# Speeding up SAT solver by exploring CNF symmetries : Revisited

Arup Kumar Ghosh
School of Electrical Engineering and Computer Science
University of Central Florida
Orlando, FL 32826
Email – akghosh@cs.ucf.edu

**Abstract**

Boolean Satisfiability solvers have gone through dramatic improvements in their performances and scalability over the last few years by considering symmetries. It has been shown that by using graph symmetries and generating symmetry breaking predicates (SBPs) it is possible to break symmetries in Conjunctive Normal Form (CNF). The SBPs cut down the search space to the non-symmetric regions of the space without affecting the satisfiability of the CNF formula. The symmetry breaking predicates are created by representing the formula as a graph, finding the graph symmetries and using some symmetry extraction mechanism (Crawford *et al.*). Here in this paper we take one non-trivial CNF and explore its symmetries. Finally, we generate the SBPs and adding it to CNF we show how it helps to prune the search tree, so that SAT solver would take short time. Here we present the pruning procedure of the search tree from scratch, starting from the CNF and its graph representation. As we explore the whole mechanism by a non-trivial example, it would be easily comprehendible. Also we have given a new idea of generating symmetry breaking predicates for breaking symmetry in CNF, not derived from Crawford's conditions. At last we propose a backtrack SAT solver with inbuilt SBP generator.

## 1   Problem description

Boolean satisfiability problem checks if a Boolean expression in CNF has some truth assignments of its variables for which the expression evaluates to true. CNF formulas are vastly used in automatic theorem proving and in Electronic Design Automation. As these formulas contain human-design artifacts and therefore often contain a great deal of symmetry, which causes satisfiability solvers to explore many redundant truth –assignments. Symmetry breaking predicates can be appended to the CNF formula to remove the symmetry while keeping satisfiability same. These Symmetry breaking predicates are created by expressing the formula as a graph, finding the graph symmetries, and then applying SBP construction algorithms on graph symmetries. The new CNF formula has the same satisfiability, but has pruned search space and therefore can be solved faster.

The problem is described by using one trivial example below: For e.g. $\theta$ is a CNF and $\theta = (X_1 + \overline{X_2})(\overline{X_1} + X_2)$.

Here $X_1, X_2, \ldots\ldots\ldots\ldots, X_n$ are variables.

$X_1, \overline{X_1}, X_2, \overline{X_2}, \ldots\ldots\ldots\ldots, X_n, \overline{X_n}$ are called literals. $(X_1 + \overline{X_2}), (\overline{X_1} + X_2)$ are clauses.

| $X_1$ | $X_2$ | $(X_1 + \overline{X_2})(\overline{X_1} + X_2)$ |
|---|---|---|
| 0 | 0 | 1 |
| 0 | 1 | 0 |
| 1 | 0 | 0 |
| 1 | 1 | 1 |

**Table 1. Truth Table of CNF $\theta$**

Now $\theta$ is called satisfiable for those values of $X_1, X_2$ for which $\theta = 1$. So, for the above case $\theta$ is satisfiable for $\{X_1 = 0, X_2 = 0\}$ and $\{X_1 = 1, X_2 = 1\}$. The required steps are shown below for finding the symmetry breaking predicates (SBP):

**Step 1:** Representation of CNF as a graph. Clauses are denoted by $C_1$ and $C_2$ where $C_1 = (X_1 + \overline{X_2})$ and $C_2 = (\overline{X_1} + X_2)$. Hence, the graph G of CNF would look like as Figure 1. The construction mechanism of graph is discussed in detail in section 4.1.

**Step 2:** Now using the graph symmetry we have to find out the symmetry breaking predicates. For finding the symmetries of the graph in terms of set of irredundant generators we have to use a suitable graph automorphism.

**Step 3:** After that we have to map the graph symmetries back to symmetries of the formula.

**Step 4:** Then we have to construct an appropriate symmetry breaking predicate (SBP) $\alpha$ and conjoint it with the formula.

**Step 5:** After that we have to solve $\theta * \alpha$ using a suitable SAT solver.

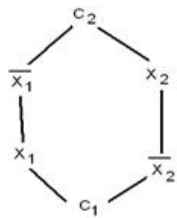

**Figure 1. Graph Form of CNF**

More detailed descriptions of this procedure are provided in consecutive sections. In Section 2 we discuss some related definitions. The above step-by-step procedure is discussed by in Section 4. Section 8 discusses about the new idea of SBP generation. Proposed SAT solver is in section 7. Section 8 and 9 include results and conclusions respectively.

## 2   Some Definitions

### 2.1   Conjunctive normal form / Variables / Literals / Clauses

In boolean logic, a formula is in Conjunctive Normal Form (CNF) if it is a conjunction of clauses, where a clause is a disjunction of literals. As a normal form, it is useful in automated theorem proving. It is similar to the canonical product of sums form used in circuit theory. All conjunctions of literals and all disjunctions of literals are in CNF, as they can be seen as conjunctions of one-literal clauses and disjunctions of a single clause, respectively. As in the Disjunctive Normal Form (DNF), the only propositional connectives a formula in CNF can contain are AND, OR, and NOT. The NOT operator can only be used as part of a literal, which means that it can only precede a propositional variable.

For example φ is a CNF and say $\varphi = (X_1 + \overline{X_2})(\overline{X_1} + X_2)$

Here $X_1, X_2, \ldots\ldots\ldots\ldots, X_n$ are variables.

$X_1, \overline{X_1}, X_2, \overline{X_2}, \ldots\ldots\ldots\ldots, X_n, \overline{X_n}$ are called literals.

$(X_1 + \overline{X_2}), (\overline{X_1} + X_2)$ are clauses.

### 2.2   Symmetries in SAT

Permutations of variables that preserve clauses create symmetries in SAT. For example, in case of CNF θ = (a + b + c)(d + e + f) if we permutate variables a, b we get the same CNF. This is shown in Table 2 along with other permutation.

| Permutation | Formula |
|---|---|
| (ab) | (b + a + c)(d + e + f) |
| (ad)(be)(cf) | (d + e + f)(a + b + c) |

**Table 2. Symmetries in CNF θ**

### 2.3   Graph Automorphism

An automorphism of a graph is a graph isomorphism with itself, i.e. a mapping from the vertices of the given graph G back to vertices of G such that the resulting graph is isomorphic with G. The sets of automorphisms define a permutation group. For every group γ, there exists a graph whose automorphism group is isomorphic to γ. The automorphism groups of a graph characterize its symmetries, and are therefore very useful in determining certain of its properties.

## 2.4 NP-complete problem

A decision problem A is NP-complete if:
  I.   A is in NP (Nondeterministic Polynomial), and
  II.  Every problem in NP is reducible to A in polynomial time.
A can be shown to be in NP by showing that a solution to the problem can be verified in polynomial time.

# 3 Symmetry Breaking Predicates Construction Algorithm

### 3.1 Lex-leader Formula [Crawford *et al.* '96]

Given a group of $\Pi = \{\Pi_1, \ldots \ldots \Pi_m\}$ defined over totally − ordered variables $X_1 < X_2 < \cdots < X_n$. For each symmetry Π, construct a permutation predicate PP(Π) in terms of bit predicate BP(Π).
According to the formula,

$$PP(\Pi) = \bigcap_{1 \leq i \leq n} BP(\Pi, i)$$

$$BP(\Pi, i) = \left[\bigcap_{1 \leq j \leq i-1} \left(x_j = x_j^\Pi\right)\right] \rightarrow (x_i \leq x_i^\Pi)$$

$$PP(\Pi) = \bigcap_{1 \leq i \leq n} \left[\bigcap_{1 \leq j \leq i-1} \left(x_j = x_j^\Pi\right)\right] \rightarrow (x_i \leq x_i^\Pi)$$

### 3.2 SBPs in terms of PP

SBPs in terms of PPs:
$$LL - SBP(\Pi) = \bigcap_{1 \leq i \leq m} PP(\Pi_i) = PP(\Pi_1).PP(\Pi_2) \ldots \ldots PP(\Pi_{m-1}).PP(\Pi_m)$$

# 4 Applying SBP construction algorithm on a CNF

For explaining the whole procedure we are considering following non-trivial CNF:
$\theta = (X_1 + \overline{X_2} + X_3 + X_4 + X_5)(X_2 + \overline{X_3} + \overline{X_4} + X_5)(\overline{X_1} + X_2 + \overline{X_5})$
There are 3 clauses in $\theta$: $C_1 = (X_1 + \overline{X_2} + X_3 + X_4 + X_5)$, $C_2 = (X_2 + \overline{X_3} + \overline{X_4} + X_5)$, $C_3 = (\overline{X_1} + X_2 + \overline{X_5})$

## 4.1 CNF to Colored Graph

Graph G is drawn such a way so that each clauses $C_1, C_2, C_3$ and literals represent vertex in graph G. Each clause representing vertex is connected by an edge to those literals which are present in that clause. Each variable is connected to their complement.

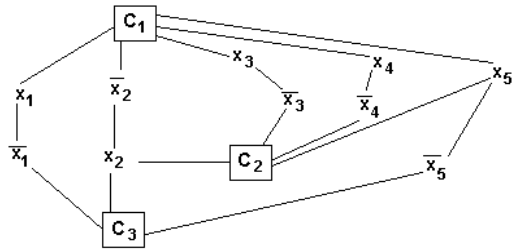

**Figure 2. Graph G = (V, E) representing the CNF $\theta$**

## 4.2 Colored Graph to CNF symmetries

Here we use NAUTY algorithm to find out the symmetry of G, $\gamma$. Let $G^\Pi$ be the relabeling of G with respect to discrete coloring $\Pi$. If $\Pi_1$ and $\Pi_2$ are discrete colorings, and $\Pi_1{}^\gamma = \Pi_2$, then $\gamma$ is symmetry of G iff $G^{\Pi_1} = G^{\Pi_2}$.

**Finding $\gamma$**
Coloring is an ordered partition of $\Pi$ of V – the cell of $\Pi$ form a sequence, not simply an unordered set.
From Figure 2 we have the initial coloring as,
$\Pi = [\, X_1\ \overline{X_1}\ X_2\ \overline{X_2}\ X_3\ \overline{X_3}\ X_4\ \overline{X_4}\ X_5\ \overline{X_5}\, |\, C_1\ C_2\ C_3\, ]$
Now the refinement steps are as follows,

$$\Pi = [\, X_1\ \overline{X_1}\ X_2\ \overline{X_2}\ X_3\ \overline{X_3}\ X_4\ \overline{X_4}\ X_5\ \overline{X_5}\, |\, C_1\ C_2\ C_3\, ]$$
$$\downarrow$$
$$[\, X_1\ \overline{X_1}\ \overline{X_2}\ X_3\ \overline{X_3}\ X_4\ \overline{X_4}\ \overline{X_5}\, |\, X_2\ X_5\, |\, C_3\, |\, C_2\, |\, C_1\, ]$$
$$\downarrow$$
$$[\, X_1\ \overline{X_1}\ X_3\ \overline{X_3}\ X_4\ \overline{X_4}\, |\, \overline{X_2}\ \overline{X_5}\, |\, X_2\ X_5\, |\, C_3\, |\, C_2\, |\, C_1\, ]$$
$$\downarrow$$
$$[\, X_1\ X_3\ \overline{X_3}\ X_4\ \overline{X_4}\, |\, \overline{X_1}\, |\, \overline{X_2}\, |\, \overline{X_5}\, |\, X_5\, |\, X_2\, |\, C_3\, |\, C_2\, |\, C_1\, ]$$

$$[X_3 \overline{X_3} \ X_4 \overline{X_4} | \ X_1 \ |\overline{X_1}\ | \overline{X_2}\ | \overline{X_5} \ | \ X_5 \ | \ X_2 \ | \ C_3 \ | \ C_2 \ | \ C_1 \ ]$$
$$\downarrow$$
$$[X_3 \ X_4 \ \overline{X_3} \ \overline{X_4} | \ X_1 \ |\overline{X_1}|\overline{X_2}\ |\overline{X_5}| \ X_5 \ | \ X_2 \ | \ C_3 \ | \ C_2| \ C_1 \ ]$$
$$\downarrow$$
$$\Pi' = [X_3 \ X_4 \ | \ \overline{X_3} \ \overline{X_4} \ | \ X_1 \ |\overline{X_1}|\overline{X_2}\ | \ X_5 \ | \ X_5 \ | \ X_2 | \ C_3 \ | \ C_2 \ | \ C_1 \ ]$$

If above refinement procedure returns a discrete coloring Π' i.e. every cell of the partition is singleton then all vertices can be distinguished, in that case the graph G must possess no symmetries besides the identity. Here $\Pi'$ is not discrete i.e. there is some non-singleton cell in Π' representing vertices that could not be distinguished based on degree – they may participate in some symmetry.

According to NAUTY algorithm it selects some non-singleton cell T of Π', called the target cell and form |T| coloring descendent from $\Pi'$, each identical to $\Pi'$ except that one t ϵ T is indivisualized in front of $T - \{t\}$. Each of these coloring is gradually refined and more descendent coloring are generated if the refined coloring are not discrete, this procedure is repeated till discrete coloring is reached.

According to above discussion we have to select some target cell. Here we select target cell $\{X_3 \ X_4\}$.

$$[X_3 \ X_4 \ | \ \overline{X_3} \ \overline{X_4}| \ X_1| \ \overline{X_2} \ \ \overline{X_5}| \ X_5 \ | \ X_2 \ | \ C_3 \ | \ C_2 \ | \ C_1 \ ] \ \overrightarrow{X_3}$$
$$[X_3 \ | \ X_4 \ | \ \overline{X_3} \ \overline{X_4}| \ X_1 \ |\overline{X_1} \ | \ \overline{X_2} \ | \ \overline{X_5}| \ X_5 \ | \ X_2 \ | \ C_3 \ | \ C_2 \ | \ C_1 \ ]$$
$$\downarrow$$
$$[X_3 \ | \ X_4 \ | \overline{X_4} \ | \ \overline{X_3}| \ X_1 \ |\overline{X_1} \ \ \overline{X_2}| \ \overline{X_5}| \ X_5 \ | \ X_2| \ C_3 \ | \ C_2 \ | \ C_1 \ ] \text{ ----------(i)}$$
$$[X_3 \ X_4| \ \overline{X_3} \ \overline{X_4}| \ X_1 \ |\overline{X_2} \ |\overline{X_5}| \ X_5 \ | \ X_2 \ | \ C_3 \ | \ C_2| \ C_1 \ ] \ \overrightarrow{X_4}$$
$$[X_4 X_3| \ \overline{X_3} \ \overline{X_4} \ | \ X_1 \ |\overline{X_1} \ | \ \overline{X_2}| \ \overline{X_5}| \ X_5 \ | \ X_2| \ C_3 \ | \ C_2 \ | \ C_1 \ ]$$
$$\downarrow$$
$$[X_4 \ | \ X_3 \ |\overline{X_3} \ | \ \overline{X_4} \ | \ X_1 \ |\overline{X_1} \ |\overline{X_2}| \ \overline{X_5}| \ X_5 \ | \ X_2 \ | \ C_3 \ | \ C_2 \ | \ C_1 \ ] \text{ ---------(ii)}$$

Using (i) and (ii) we have , $\gamma = (X_3 \ X_4)(\overline{X_3} \ \overline{X_4})$

## 4.3 Forming Symmetry Breaking Predicates (SBPs)

Using Lex-leader Formula formula in our non-trivial example we have,
$\gamma = (X_3 \ X_4)(\overline{X_3} \ \overline{X_4})$
$\Pi = \begin{Bmatrix} X_1 & X_2 & X_3 & X_4 & X_5 \\ X_1 & X_2 & X_4 & X_3 & X_5 \end{Bmatrix}$
Here, Π is the permutation of the set of 2n literals L = $\{X_1, \overline{X_1}, \dots \dots, X_n, \overline{X_n}\}$ is a function Π : L → L that is both one-to-one and onto. If $X_j$ is the image of $X_i$ under Π, then it is denoted by $X_j = X_i^{\Pi}$.

When $\overline{X_i}$ maps to $\overline{X_j}$ that equivalent to mapping from $X_i$ to $X_j$ due to preservation of Bolean consistency.

$BP(\Pi, 1) = X_1 \leq X_1 = 1$ //maps to corresponding bit itself
$BP(\Pi, 2) = (X_1 = X_1) \rightarrow (X_2 \leq X_2) = 1$ //maps to corresponding bit itself
$BP(\Pi, 3) = (X_1 = X_1)(X_2 = X_2) \rightarrow (X_3 \leq X_4) = \overline{X_3} + X_4$
$BP(\Pi, 4) = (X_1 = X_1)(X_2 = X_2)(X_3 = X_4) \rightarrow (X_4 \leq X_3) = 1$ //cycle ends
$BP(\Pi, 5) = (X_1 = X_1)(X_2 = X_2)(X_3 = X_4)(X_4 = X_3) \rightarrow (X_5 \leq X_5) = 1$
//maps to corresponding bit itself

$$\text{Hence, PP}(\Pi) = \bigcap_{1 \leq i \leq n} BP(\Pi, i) = \overline{X_3} + X_4$$

### 4.4 CNF with SBPs

CNF with SBP $(\theta') = \theta * SBP = \theta = (X_1 + \overline{X_2} + X_3 + X_4 + X_5)(X_2 + \overline{X_3} + \overline{X_4} + X_5)(\overline{X_1} + X_2 + \overline{X_5}) \; (\overline{X_3} + X_4)$

## 5 Pruning Using SBP

The truth table for above CNF is given in Table 3. This table is showing different combination of boolean variables and corresponding values of the CNF and values of SBPs. Here we can see that for some cases SBPs are equal to 0 and other cases it is equal to 1. When the value = 0, that means some row of the truth table has the symmetric assignment. When the modified CNF is given to SAT solver it will discard those assignments for which SBP = 0. In the last column of the table the symmetric assignments are indexed by same numbers.

For finding out the satisfiability of a CNF, backtrack SAT solver goes through the search tree. When symmetry is there between the assignments then for testing the satisfiability, one among all symmetric assignments is sufficient, which saves the searching time and speed up the overall calculation.

| $X_1$ | $X_2$ | $X_3$ | $X_4$ | $X_5$ | $C_1$ | $C_2$ | $C_3$ | $\theta$ | SBP | |
|---|---|---|---|---|---|---|---|---|---|---|
| 0 | 0 | 0 | 0 | 0 | 1 | 1 | 1 | 1 | 1 | |
| 0 | 0 | 0 | 0 | 1 | 1 | 1 | 1 | 1 | 1 | |
| 0 | 0 | 0 | 1 | 0 | 1 | 1 | 1 | 1 | 1 | ① |
| 0 | 0 | 0 | 1 | 1 | 1 | 1 | 1 | 1 | 1 | ② |
| 0 | 0 | 1 | 0 | 0 | 1 | 1 | 1 | 1 | 0 | ① |
| 0 | 0 | 1 | 0 | 1 | 1 | 1 | 1 | 1 | 0 | ② |
| 0 | 0 | 1 | 1 | 0 | 1 | 0 | 1 | 0 | 1 | |
| 0 | 0 | 1 | 1 | 1 | 1 | 1 | 1 | 1 | 1 | |
| 0 | 1 | 0 | 0 | 0 | 0 | 1 | 1 | 0 | 1 | |
| 0 | 1 | 0 | 0 | 1 | 1 | 1 | 1 | 1 | 1 | |
| 0 | 1 | 0 | 1 | 0 | 1 | 1 | 1 | 1 | 1 | ③ |
| 0 | 1 | 0 | 1 | 1 | 1 | 1 | 1 | 1 | 1 | ④ |
| 0 | 1 | 1 | 0 | 0 | 1 | 1 | 1 | 1 | 0 | ③ |
| 0 | 1 | 1 | 0 | 1 | 1 | 1 | 1 | 1 | 0 | ④ |
| 0 | 1 | 1 | 1 | 0 | 1 | 1 | 1 | 1 | 1 | |
| 0 | 1 | 1 | 1 | 1 | 1 | 1 | 1 | 1 | 1 | |
| 1 | 0 | 0 | 0 | 0 | 1 | 1 | 1 | 1 | 1 | |
| 1 | 0 | 0 | 0 | 1 | 1 | 1 | 0 | 0 | 1 | |
| 1 | 0 | 0 | 1 | 0 | 1 | 1 | 1 | 1 | 1 | ⑤ |
| 1 | 0 | 0 | 1 | 1 | 1 | 1 | 0 | 0 | 1 | ⑥ |
| 1 | 0 | 1 | 0 | 0 | 1 | 1 | 1 | 1 | 0 | ⑤ |
| 1 | 0 | 1 | 0 | 1 | 1 | 1 | 0 | 0 | 0 | ⑥ |
| 1 | 0 | 1 | 1 | 0 | 1 | 0 | 1 | 0 | 1 | |
| 1 | 0 | 1 | 1 | 1 | 1 | 1 | 0 | 0 | 1 | |
| 1 | 1 | 0 | 0 | 0 | 1 | 1 | 1 | 1 | 1 | |
| 1 | 1 | 0 | 0 | 1 | 1 | 1 | 1 | 1 | 1 | |
| 1 | 1 | 0 | 1 | 0 | 1 | 1 | 1 | 1 | 1 | ⑦ |
| 1 | 1 | 0 | 1 | 1 | 1 | 1 | 1 | 1 | 1 | ⑧ |
| 1 | 1 | 1 | 0 | 0 | 1 | 1 | 1 | 1 | 0 | ⑦ |
| 1 | 1 | 1 | 0 | 1 | 1 | 1 | 1 | 1 | 0 | ⑧ |
| 1 | 1 | 1 | 1 | 0 | 1 | 1 | 1 | 1 | 1 | |
| 1 | 1 | 1 | 1 | 1 | 1 | 1 | 1 | 1 | 1 | |

**Table 3. Truth Table of non-trivial** CNF $\theta = (X_1 + \overline{X_2} + X_3 + X_4 + X_5)(X_2 + \overline{X_3} + \overline{X_4} + X_5)(\overline{X_1} + X_2 + \overline{X_5})$

**Figure 3. Search Tree For the CNF. Here round cornered boxes indicate those parts of the tree which are discarded by the SBP.**

In Figure 3, the search tree for $\theta = (X_1 + \overline{X_2} + X_3 + X_4 + X_5)(X_2 + \overline{X_3} + \overline{X_4} + X_5)(\overline{X_1} + X_2 + \overline{X_5})$ is drawn. Without adding the SBP, SAT solver goes through the entire tree to find out satisfiable assignments. But when we add SBP with it it discard those part of the tree for which SBP value equals to zero. For our example the discarded parts of the tree are shown by a box with rounded corners. As the volume of the search tree is decreased, the SAT would take less time than searching through original search tree.

## 6 Another idea of generation of SBP

| $X_1$ | $X_2$ | $X_3$ | $X_4$ | $X_5$ | $C_1$ | $C_2$ | $C_3$ | $\theta$ | SBP | |
|---|---|---|---|---|---|---|---|---|---|---|
| 0 | 0 | 0 | 0 | 0 | 1 | 1 | 1 | 1 | 1 | |
| 0 | 0 | 0 | 0 | 1 | 1 | 1 | 1 | 1 | 1 | |
| 0 | 0 | 0 | 1 | 0 | 1 | 1 | 1 | 1 | 1 | ① |
| 0 | 0 | 0 | 1 | 1 | 1 | 1 | 1 | 1 | 1 | ② |
| 0 | 0 | 1 | 0 | 0 | 1 | 1 | 1 | 1 | 0 | ③ |
| 0 | 0 | 1 | 0 | 1 | 1 | 1 | 1 | 1 | 0 | ④ |
| 0 | 0 | 1 | 1 | 0 | 1 | 0 | 1 | 0 | 1 | |
| 0 | 0 | 1 | 1 | 1 | 1 | 1 | 1 | 1 | 1 | |

**Table 4. Part of table 3**

Here in Table 4 we have some part of the truth table. Now SBP eliminates redundant assignments, by making SBP value = 0. In the previous work SBP becomes 0 for assignments in long black box in Table 4. Without using Lex-leader formula we can make our SBP, so that it can be 0 for $X_3 = 1$ and $X_4 = 0$ and for other cases SBP is non-zero.
Now for two variables (a, b), truth table for $\bar{a} + b$ is given below:

| a | b | $\bar{a} + b$ |
|---|---|---|
| **0** | 0 | 1 |
| 0 | 1 | 1 |
| 1 | 0 | 0 |
| **1** | 1 | 1 |

**Table 5. Truth table for $\bar{a} + b$**

We are interested about gray portion of the truth table.
If we relate our motto to the above truth table, then we need to add $\overline{X_3} + X_4$ as the SBP and its value would become 0 for $X_3 = 1$ and $X_4 = 0$.
Hence when, $\Pi = (X_3 X_4)(\overline{X_3}\ \overline{X_4})$
SBP becomes $\overline{X_3} + X_4$.
In general when $\Pi = (X_1 X_2)(X_3 X_4)(X_5 X_6)$, SBP becomes $(\overline{X_1} + X_2)(\overline{X_3} + X_4)(\overline{X_5} + X_6)$.

## 6.1 Algorithm to generate SBP manually

```
SBP_Partial = 1; SBP_Total = 1;
for every Π
{
            for every mapped pair nodes (Say X_i and X_j) in Π
                    {
                    SBP_Partial = SBP_Partial * ( $\overline{X_i}$ + $X_j$ );
                    }
            //SBP_Total is the Final SBP
            SBP_Total = SBP_Total * SBP_Partial;
}
```

This approach needs less computation than the previous approach. So, this approach should be faster.

## 6.2 Applying above algorithm on an example

Say, $\Pi = ( X_1 \ X_2 ) ( X_3 \ X_4 ) ( X_5 \ X_6 )$
Here we have one Π. Hence first "for loop" will iterate 1 time.
In Π we have 3 mapped pair nodes. Hence next will iterate 3 times
Now using above algorithm we have,
1 SBP_Partial = 1; SBP_Total = 1;
2 SBP_Partial = 1 * ( $\overline{X_1}$ + $X_2$ ) = ( $\overline{X_1}$ + $X_2$ );
3 SBP_Partial = ( $\overline{X_1}$ + $X_2$ )( $\overline{X_2}$ + $X_3$ )
4 SBP_Partial = ( $\overline{X_1}$ + $X_2$ )( $\overline{X_2}$ + $X_3$ )( $\overline{X_4}$ + $X_5$ )
5 SBP_Total = SBP_Total * SBP_Partial = ( $\overline{X_1}$ + $X_2$ )( $\overline{X_2}$ + $X_3$ )( $\overline{X_4}$ + $X_5$ )

# 7 Proposed SAT solver

In the figure 4 we have shown the different steps of SBPs generation and the way to find out solution using backtrack search. As there might have lot of symmetries in CNF with thousands of variables and clauses, it will be very costly to find out SBPs using Crawford et al. If proposed manual method for SBP generation is used, then it will be less costlier. Hence, the total time needed to solve SAT problem starting from original CNF formula, considering the SBPs generation, will be reduced and hence improve the overall performance.
Again proposed backtrack SAT solver with inbuilt SBP generator will reduce work and time. In that case we don't have to think about various part of overall process, described in the following figure. We just need to provide the original CNF to the system and the system will take care of everything, starting from graph representation to finding out the solution.

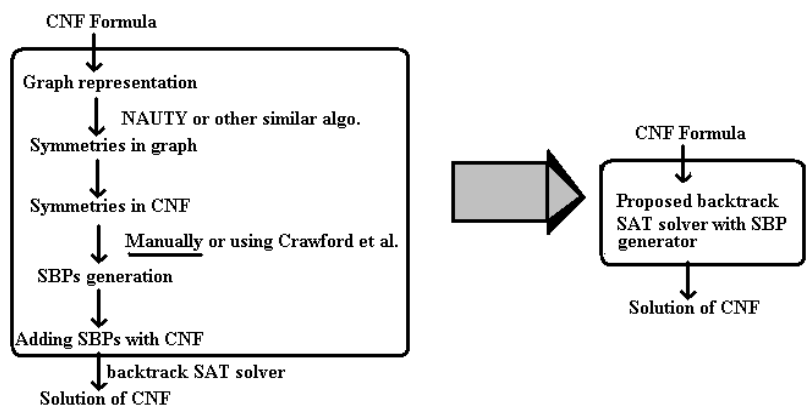

**Figure 4. Flow diagram of solving SAT problem using SBPs. Proposed manual construction of symmetry breaking predicates improves the performance than others known methods, as this method does not need any complicated calculation. On the right hand side of above diagram the proposed backtrack SAT solver with inbuilt SBP generator is shown.**

## 8 Results

The example in section 7 deals with 5 variables CNF. Using symmetry breaking 8 symmetric assignments are found. Without symmetry detection roughly $2^5 = 32$ assignments have to be examined at worst case to test the satisfiability. But using Symmetry it drops down to 24. The effect of this reduction is clearly shown in Figure 3. For CNF with large no. of variables Symmetry Breaking Predicates help to prune the search tree significantly. It is proved that solving satisfiability is NP complete, i.e. no known polynomial time algorithm is there. In this scenario somehow decreasing the search space and improving the run time is good idea.

Again proposed manual SBP generation improves the overall performance by reducing the computation.

## 9 Conclusions and future work

This paper shows how symmetry in CNF can be broken and adding SBPs with CNF how it speeds up the SAT solver. In this paper every step has been done manually, so that it can be clearly understandable. There are few paper available on this topic but no paper has described the whole story in one place. Here we have tried to make it one place and tried to use one non trivial example so that the main idea behind this method can be revealed.

We already proposed the manual SBP generation and a backtrack SAT solver with inbuilt SBP generators.
Our future work includes full implementation of this backtrack SAT solver with inbuilt SBP generators.